\documentclass[12pt, twoside, leqno]{article}
\usepackage{amsmath,amsthm,amscd, amsfonts, amssymb, graphicx, color}
\usepackage{amssymb,latexsym}
\usepackage{enumerate}
\usepackage{hyperref}
\newtheorem{theorem}{Theorem}[section]
\newtheorem{lemma}[theorem]{Lemma}

\newtheorem{corollary}[theorem]{Corollary}
\theoremstyle{definition}
\newtheorem{definition}[theorem]{Definition}
\newtheorem{example}{Example}

\newtheorem{question}[theorem]{Question}

\numberwithin{equation}{section}
\begin{document}
\baselineskip=17pt

\title{On a question related to bounded approximate identities of ideals in Banach algebras}

\author{Mohammad Fozouni\\
\small{{Department of Mathematics, Faculty of Sciences and Engineering,}}\\
\small{{P. O. Box 163, Gonbad Kavous, Golestan, Iran,}}\\
\small{{E-mail: fozouni@gonbad.ac.ir or fozouni@hotmail.com}}}
\date{}
\maketitle

\begin{abstract}
In this paper we give an example of a Banach algebra $A$ and a closed ideal $I$ of $A$ such that the multiplier algebra of $I$ is equal to $A$ but $I$ does not have any bounded approximate identity. In the case that $I$ has an approximate identity, we give a necessary condition on $I$ for which $A=\mathcal{M}(I)$, where $\mathcal{M}(I)$ denotes the multiplier algebra of $I$. Finally, as a corollary of our results, we show that the Fourier algebra of an amenable group is strictly dense in the Fourier-Stieltjes algebra.

\vspace{20pt}
\noindent\textbf{MSC 2010:} {46H05, 22D15}\\
\noindent\textbf{Keywords:} {Banach algebra, approximate identity,multiplier, Fourier algebra}
\end{abstract}
\section{Introduction}
This work motivated by the following observation:

We know that the multiplier algebra of the group algebra is the measure algebra and the multiplier algebra of the Fourier algebra is the Fourier-Stieltjes algebra when the underlying group is amenable. Clearly, the group algebra and the Fourier algebra of an amenable locally compact group both have bounded approximate identity. So, we ask the following questions:

If $A$ is a Banach algebra and $I$ is a closed ideal of $A$ such that the multiplier algebra of $I$ is $A$, that is $\mathcal{M}(I)=A$, is it true that $I$ has a bounded approximate identity? Is the ideal $I$ unique in the representation $\mathcal{M}(I)=A$?

In the next section we try to give answer to these questions. Also, under some conditions, we give a necessary and sufficient condition on ideal $I$ of $A$ such that $\mathcal{M}(I)=A$.

\section{Main results}
Suppose that $A$ is a Banach algebra. We say that the bounded linear operator $T: A\rightarrow A$ is a (left) multiplier of $A$ if $T(ab)=T(a)b$ for each $a, b\in A$. Let $\mathcal{M}(A)$ denote the set of all multipliers of $A$.

A net $\{e_\alpha\}$ in $A$ is called  an approximate identity (a.i.) for $A$ if for every $a\in A$, $\|ae_\alpha-a\|+\|e_\alpha a-a\|\rightarrow 0$. If the net $\{e_\alpha\}$ is bounded and has the mentioned property we say that $\{e_\alpha\}$ is a bounded approximate identity (b.a.i.) for $A$.

To see the standard definitions of the undefined concepts in the sequel one can see the references \cite{Dales,  Murphy, Pier}.

Now, we give the following definition for the sake of convenience in our notations.
\begin{definition}
Let $A$ be a Banach algebra. We say that $A$ is a multiplierly generated algebra (briefly MGA) if there exists a closed ideal $I$ in $A$ such that $A=\mathcal{M}(I)$, that is, the mapping $a\longrightarrow L_a$ from $A$ to $\mathcal{M}(I)$ is an isometric isomorphism, where $L_a(b)=ab$.
\end{definition}
One can see easily that if $A$ is a unital algebra, then $\mathcal{M}(A)=A$. So, we are interested in the case that $I\subsetneq A$ and $\mathcal{M}(I)=A$.

\begin{example}\label{7}
Let $G$ be a locally compact group. Then $M(G)$ is a MGA by Wendel's Theorem. More precisely, $\mathcal{M}(L^1(G))=M(G)$. To see the definitions of the measure algebra $M(G)$ and the group algebra $L^1(G)$ we refer the reader to \cite[Section 3.3]{Dales}.
\end{example}
Suppose that $G$ is a locally compact group and suppose that $L^\infty (G)$ is the space of all essentially bounded and Borel measurable functions from $G$ into $\mathbb{C}$. The group $G$ is said to be amenable if there exists an $m\in L^{\infty}(G)^{*}$ such that $m\geq 0$, $m(1)=1$ and $m(L_{x}f)=m(f)$ for each $x\in G$ and $f\in L^{\infty}(G)$ where $L_{x}f(y)=f(x^{-1}y)$. Suppose that $A(G)$ denotes the Fourier algebra and $B(G)$ denotes the Fourier-Stieltjes algebra of $G$. We know that $A(G)$ is a closed ideal of $B(G)$; see \cite{Pier} for more details of the amenable groups and the Fourier algebra.
\begin{example}\label{3}
 Let $G$ be a locally compact group. Then $G$ is amenable if and only if $\mathcal{M}(A(G))=B(G)$ , so $B(G)$ is a MGA; see \cite[Theorem 1]{Losert}.  
\end{example}
\begin{example}
  Let $X$ be a locally compact Hausdorff space. Then by \cite[Example 1.4.13]{Kaniuth2} we know that $\mathcal{M}(C_0(X))=C_b(X)$. So, $C_b(X)$ is a MGA.
\end{example}
\begin{example}\label{6}
  Let $H$ be a Hilbert space and $K(H)$ be the space of compact operators. Using the definition of compact operators, one can see that $K(H)$ is a closed ideal of $B(H)$. Also, we know that $\mathcal{M}(K(H))=B(H)$; see \cite[Example 3.1.2]{Murphy}, note that in this example double centralizer obtained, but similarly one can prove the result for  multipliers. Therefore, $B(H)$ is a MGA.
\end{example}
\begin{example}
Every non-unital Banach algebra $A$ is not a MGA, because $\mathcal{M}(I)$ is unital and this forces $A$ to be unital.
\end{example}

Examples \ref{7} to \ref{6} show that for a large number of important Banach algebras $A$  and its closed ideals $I$, if $\mathcal{M}(I)=A$, then $I$ has a b.a.i. Now, we ask the following questions and try to answer these questions in the sequel.
\begin{question}\label{1}
If $A$ is a MGA with $\mathcal{M}(I)=A$, is it true that $I$ should has a  b.a.i.?
\end{question}

\begin{question}\label{2}
 If $\mathcal{M}(I)=\mathcal{M}(J)=A$ where $I, J$ are two proper closed ideals of $A$, is $I=J$ or it is not necessarily true?
\end{question}

The following example shows that there exists a Banach algebra $A$ and an ideal $I$ of $A$ such that $\mathcal{M}(I)=A$, .i.e., $A$ is multiplierly generated but $I$ does not have any b.a.i. But Question \ref{1} does not fail, because one can easily check that $I$ is not a closed ideal of $A$.

\begin{example}
Let $A=\ell^1(\mathbb{N})$ and consider $\ell^{\infty}(\mathbb{N})$ with the sup-norm, i.e., for $f\in \ell^{\infty}(\mathbb{N})$, $\|f\|_{\infty}=\sup_{n\in \mathbb{N}} |f(n)|$. One can see that $\mathcal{M}(A)=\ell^{\infty}(\mathbb{N})$, that is, $\mathcal{M}(A)$ is isometrically isomorphic to $\ell^{\infty}(\mathbb{N})$ where $g\in \ell^{\infty}(\mathbb{N})$ acts on $A$ by pointwise multiplication; see \cite[Exercise 1.6.43]{Kaniuth2}. Indeed, the mapping $g\rightarrow L_g$ is an isometric isomorphism between $\ell^{\infty}(\mathbb{N})$ and $\mathcal{M}(A)$. To see this, let $T\in \mathcal{M}(A)$. It is easily verified that for each $n\in \mathbb{N}$, there exits $f\in A$ such that $f(n)\geq 1$ and $\|f\|_1\leq 2$. Now, put $g(n)=\frac{(Tf)(n)}{f(n)}$ where $f$ is in  $A$, $f(n)\geq 1$ and $\|f\|_1\leq 2$. Since $T$ is a multiplier, $g$ is a well-defined function and it is not depend on $f$.
Also,  $g\in \ell^\infty(\mathbb{N})$ and $T=L_g$, since if $f(n)=0$, then $(Tf)(n)^2=f(n)(Tf)(n)=0$ and hence $(Tf)(n)=0$. On the other hand, for each $n\in \mathbb{N}$, if $f(n)=1$ and $\|f\|_1=1$, then we have

$|g(n)|\leq |g(n)f(n)|=|Tf(n)|\leq \|Tf\|\leq\|T\|.$
So, $\|T\|=\|g\|_\infty$ and this shows that the mapping is an isometry.
 Also, one can see easily that $A$ is an ideal of $\ell^{\infty}(\mathbb{N})$ but it is not closed and $A$ does not have any b.a.i., it has only an a.i. .
\end{example}
In the sequel we give a negative answer to Questions \ref{1} and \ref{2}. But first we give the following lemma.
\begin{lemma}\label{5}
  Suppose that $I\subseteq J\subseteq A$ where $I, J$ are closed ideals of Banach algebra $A$, $I$ has an approximate identity and $A$ is unital. Then
  $$\mathcal{M}(A)\subseteq \mathcal{M}(J)\subseteq \mathcal{M}(I).$$
\end{lemma}
\begin{proof}
  Let $T\in \mathcal{M}(A)$. We show that $T_{|J}$ is a function from $J$ into $J$ and this shows that $T_{|J}\in \mathcal{M}(J)$. Since $A$ is unital we have $\mathcal{M}(A)=A$, therefore, for each $T\in \mathcal{M}(A)$, there exists $a\in A$ such that $T=L_a$. So $T_{|J}=(L_{a})_{|J}: J\rightarrow A$ and for all $j\in J$ we have
  $$T_{|J}(j)=(L_a)_{|J}(j)=aj\in J.$$
  Hence, $T_{|J}(J)\subseteq J$.

  Now, let $\{e_\alpha\}$ be an approximate identity for $I$ and $S\in \mathcal{M}(J)$. So, $S_{|I}: I\rightarrow J$ has a  multiplier property, i.e., for $i_1, i_2\in I$, $S_{|I}(i_1i_2)=S_{|I}(i_1)i_2$
     and for each $i\in I$ we have
  $$S_{|I}(i)=\lim_\alpha S_{|I}(e_\alpha i)=\lim_\alpha e_\alpha S_{|I}(i)\in I.$$
  Hence $S_{|I}(I)\subseteq I$ and this completes the proof.
\end{proof}
Let $G$ be a locally compact group and $M(G)$ be the measure algebra consisting of all the complex-valued regular Borel measures on $G$. We have the following subspaces of $M(G)$:
\begin{align*}
  M_c(G)= & \text{ the space of all the continuous measures,}\\
    M_a(G)= & \text{ absolutely continuous measures
    respect to the Haar measure,}\\
  M_{cs}(G)= & \text{  measures which are singular respect to the Haar measure.}
\end{align*}
As in \cite[Section 3.3]{Dales},  $M_c(G)=M_a(G)\oplus_1 M_{cs}(G)$. So, $M_a(G)\subseteq M_c(G)\subseteq M(G)$. Also, $M_c(G)$ is a closed ideal of $M(G)$ and $L^1(G)$ is isometrically isomorphic to $M_a(G)$.

Now, using $M(G)$ and Lemma \ref{5} we give answer to Questions \ref{1} and \ref{2}.\\
Since $M_a(G)\subseteq M_c(G)\subseteq M(G)$, $M_a(G), M_c(G)$ are closed ideals of $M(G)$  and $M_a(G)$ has a b.a.i., by Lemma \ref{5} $\mathcal{M}(M(G))\subseteq \mathcal{M}(M_c(G)) \subseteq \mathcal{M}(M_a(G))$. Therefore, by Wendel's theorem we have $\mathcal{M}(M_c(G))=\mathcal{M}(M_a(G))=M(G)$. But by \cite[Theorem 2.7]{DGH}, if $G$ is a non-discrete locally compact group (for example $G=\mathbb{R}$, the additive group of real numbers is  non-discrete), then $\overline{M_c(G)^2}$ has infinite codimension in $M_c(G)$, that is, $\dim(\frac{M_c(G)}{\overline{M_c(G)^2}})$ is not finite and so $ M_c(G)\neq \overline{M_c(G)^2}$. Therefore, $M_c(G)$ does not have any approximate identity. Hence, Question \ref{1} and \ref{2} both fails, i.e., when $G$ is non-discrete, there exists two proper closed ideals $I, J$ of $M(G)$ such that $\mathcal{M}(I)=\mathcal{M}(J)=M(G)$ and $J$ does not have any b.a.i. .

In the sequel we present a necessary and sufficient condition for a MGA and using this assertion we show that $A(G)$ is strictly dense in $B(G)$ if $G$ is amenable.

Recall that if $\mathcal{A}$ is a Banach algebra which is contained as a closed ideal in Banach algebra $\mathcal{B}$ then the strict topology on $\mathcal{B}$ from $\mathcal{A}$ is the topology generated by the following family of semi-norms:
$$p_a(b)=\|ab\|+\|ba\| \ \ \  (b\in \mathcal{B}, a\in \mathcal{A}).$$

So, a net $\{b_\alpha\}$ in $\mathcal{B}$ tends to $b_0$ in the strict topology if for all $a\in \mathcal{A}$, we have $\|(b_\alpha-b_0)a\|, \|a(b_\alpha-b_0)\|\rightarrow 0$. We denote the closure of $\mathcal{A}$ in the strict topology by $\overline{\mathcal{A}}^{s.t.}$. We say that $\mathcal{A}$ is boundedly strictly dense in $\mathcal{B}$ if for each $b\in \mathcal{B}$ there exists a bounded net $\{e_\alpha\}$ in $A$ such that $e_\alpha\rightarrow b$ strictly. We denote this type of density by $\overline{\mathcal{A}}^{b.s.t.}=\mathcal{B}$.

\begin{theorem}\label{4}
  If $A$ is a commutative MGA with $\mathcal{M}(I)=A$ and $I$ has an approximate identity, then $I$ is strictly dense in $A$. Conversely, if $\overline{I}^{b.s.t.}=A$, then $\mathcal{M}(I)=A$.
\end{theorem}
\begin{proof}
In view of the hypothesis $\mathcal{M}(I)=A$, suppose that $K: A\rightarrow \mathcal{M}(I)$ defined by $K(a)=L_a$ is an isometric isomorphism. Hence $K$ has an inverse function $F: \mathcal{M}(I)\rightarrow A$ such that $F$ is bijective and $K\circ F=id_{\mathcal{M}(I)}$.

Now, let $a\in A$ be arbitrary. So, there exists $T\in \mathcal{M}(I)$ such that $F(T)=a$. Hence
$$T=K\circ F(T)=K(a)=L_a.$$
If $\{e_\alpha\}$ is an approximate identity for $I$ and $i\in I$, then we have
$$ai=L_a(i)=T(i)=\lim_\alpha T(e_\alpha i)=\lim_\alpha T(e_\alpha)i.$$
Therefore, $\lim_\alpha \|(a-T(e_\alpha))i\|=\lim_\alpha \|i(a-T(e_\alpha))\|=0$ for all $i\in I$ and this implies that $a\in \overline{I}^{s.t.}$. So, $I$ is strictly dense in $A$.

To prove the converse, let $T\in \mathcal{M}(I)$. We extend $T$ to a multiplier of $A$ as follows:

Define $\widetilde{T} : A\rightarrow A$ by $\widetilde{T}(a)=s.t.-\lim_\alpha T(e_\alpha)$, i.e., $\widetilde{T}(a)$ is the limit of the net $\{T(e_\alpha)\}$ in the strict topology and $\{e_\alpha\}$ is a bounded net in $I$ such that strictly tends to $a$. Clearly, $\widetilde{T}$ is a well-defined function.
Suppose that $a, b\in A$ and $\{e_\alpha\}, \{e^{'}_\beta\}$ are two nets for which tend to $a$ and $b$, respectively in the strict topology. Now, we have
\begin{align*}
  \widetilde{T}(ab) & =s.t.-\lim_{(\alpha, \beta)}T(e_\alpha e^{'}_\beta) \\
  & =s.t.-\lim_{(\alpha, \beta)}T(e_\alpha)e^{'}_\beta \\
   & =\left(s.t.-\lim_{\alpha}T(e_\alpha)\right)\left(s.t.-\lim_{\beta}e^{'}_\beta\right)\\
   & = \widetilde{T}(a)b.
\end{align*}
Note that the boundedness of $\{e_\alpha\}$ yields the third equation in the above calculations. Hence, $\widetilde{T}$ is in $\mathcal{M}(A)=A$ and so there exists $a\in A$ for which $\widetilde{T}=L_a$. Therefore, $T=\widetilde{T}_{|I}=L_a$ and this shows that $A=\mathcal{M}(I)$.
\end{proof}
In Example \ref{3} we see that $\mathcal{M}(A(G))=B(G)$ when $G$ is an amenable group. So, by Theorem \ref{4} we conclude the following corollary.
\begin{corollary}
If $G$ is an amenable group, then $\overline{A(G)}^{s.t.}=B(G)$.
\end{corollary}






\bibliographystyle{plain}
\bibliography{References}
\end{document}